%AMS TeX file for the paper
%
%     Spectral analysis and the Haar functional
%          on the quantum $SU(2)$ group
%
%by H.T.~Koelink and J.~Verding
%%%%%%%%%%%%%%%AMS TeXfile%%%%%%%%%%%%%%%%%%%%%%%%%%%%%%%%%%%%%%%%%%%%
\input amstex.tex
\documentstyle{amsppt}
\magnification=1200
\baselineskip=13pt
\hsize=6.5truein
\vsize=8.9truein
%%%%%%%%%%%%%%%%% M a c r o s %%%%%%%%%%%%%%%%%%%%%%%%%%%%%%%%%%%%%%%%%%%
%Section, Theorem and Formula Numbering%%%%%%%%%%%%%%%%%%%%%%%%%%%%%%%%%%
%%%%%%%%%%%%%%%%%%%%%%%%%%%%%%%%%%%%%%%%%%%%%%%%%%%%%%%%%%%%%%%%%%%%%%%%%
\countdef\sectionno=1
\countdef\eqnumber=10
\countdef\theoremno=11
\countdef\countrefno=12
\countdef\cntsubsecno=13
\sectionno=0
\def\newsection{\global\advance\sectionno by 1
                \global\eqnumber=1
                \global\theoremno=1
                \global\cntsubsecno=0
                \the\sectionno}

\def\newsubsection#1{\global\advance\cntsubsecno by 1
                     \xdef#1{{\S\the\sectionno.\the\cntsubsecno}}
                     \ \the\sectionno.\the\cntsubsecno.}

\def\theoremname#1{\the\sectionno.\the\theoremno
                   \xdef#1{{\the\sectionno.\the\theoremno}}
                   \global\advance\theoremno by 1}

\def\eqname#1{\the\sectionno.\the\eqnumber
              \xdef#1{{\the\sectionno.\the\eqnumber}}
              \global\advance\eqnumber by 1}

\def\thmref#1{#1}

\global\countrefno=1

\def\refno#1{\xdef#1{{\the\countrefno}}\global\advance\countrefno by 1}

%%%%%%%%%%%%%%%%%%%%%%%%%%%%Abbreviations%%%%%%%%%%%%%%%%%%%%%%%%%%%%%%%%
\def\R{{\Bbb R}}
\def\N{{\Bbb N}}

\def\Zp{{\Bbb Z}_+}
\def\hf{{1\over 2}}

\def\A{A_q(SU(2))}
\def\Hi{\ell^2(\Zp)}
\def\a{\alpha}
\def\b{\beta}
\def\g{\gamma}
\def\d{\delta}
\def\s{\sigma}
\def\t{\tau}
\def\l{\lambda}
\def\m{\mu}

\def\r{\rho_{\t ,\s}}

\def\rti{\rho_{\t ,\infty}}
\def\p{\pi}
\def\vp{\varphi}
\def\L{\Cal L}
\def\M{\Cal M}
%%%%%%%%%%%%%%%%%%%%%%%%%Reference Numbering%%%%%%%%%%%%%%%%%%%%%%%%%
\refno{\AskeI}
\refno{\AskeIMAMS}
\refno{\AskeRS}
\refno{\AskeW}
\refno{\Bere}
\refno{\Bres}
\refno{\Chih}
\refno{\Domb}
\refno{\DunfS}
\refno{\GaspR}
\refno{\KoelAF}
\refno{\KoelAAM}
\refno{\Koelpp}
\refno{\KoorIM}
\refno{\KoorOPTA}
\refno{\KoorZSE}
\refno{\MasuMNNU}
\refno{\Noum}
\refno{\NoumMCM}
\refno{\VaksS}
\refno{\WoroRIMS}
\refno{\WoroCMP}

%%%%%Beginning of the text%%%%%%%%%%%%%%%%%%%%%%%%%%%%%%%%%%%%%%%%%%%%
\topmatter
\title Spectral analysis and the Haar functional on the quantum $SU(2)$
group\endtitle
\rightheadtext{Spectral analysis and the Haar functional}
\author H.T.~Koelink and J.~Verding\endauthor
%\affil  Katholieke Universiteit Leuven\endaffil
\address Departement Wiskunde, Katholieke Universiteit Leuven,
Celestijnenlaan 200 B, B-3001 Leuven (Heverlee), Belgium\endaddress
\email Erik.Koelink\@wis.KULeuven.ac.be,
Jan.Verding\@wis.KULeuven.ac.be\endemail
\date December \the\day, 1994\enddate
\thanks The first author is supported by a Fellowship of the
Research Council of the Katholieke Universiteit Leuven. \endthanks
\abstract
The Haar functional on the quantum $SU(2)$ group is the analogue of
invariant integration on the group $SU(2)$. If restricted to a
subalgebra generated by a self-adjoint element the Haar functional can
be expressed as an integral with a continuous measure or with a discrete
measure or by a combination of both.
These results by Woronowicz and Koornwinder
have been proved by using the corepresentation
theory of the quantum $SU(2)$ group and Schur's orthogonality relations
for matrix elements of irreducible unitary corepresentations.
These results are proved here by
using a spectral analysis of the generator of the subalgebra. The
spectral measures can be described in terms of the orthogonality
measures of orthogonal polynomials by using the theory of Jacobi
matrices.
\endabstract
\endtopmatter
\document

%%%%%%%%%%%%%%%%%%%%%%%%%%%%%%%%%%%%%%%%%%%%%%%%%%%%%%%%%%%%%%%%%%%%%%
% NEW SECTION %%%%%%%%%%%%%%%%%%%%%%%%%%%%%%%%%%%%%%%%%%%%%%%%%%%%%%%%
%%%%%%%%%%%%%%%%%%%%%%%%%%%%%%%%%%%%%%%%%%%%%%%%%%%%%%%%%%%%%%%%%%%%%%

\head\newsection . Introduction \endhead

The existence of the Haar measure for locally compact groups is a
corner stone in harmonic analysis. The situation for general quantum
groups is not (yet) so nice, but for compact matrix quantum groups
Woronowicz \cite{\WoroCMP, thm.~4.2} has proved that a suitable
analogue of the Haar measure exists. This analogue of the Haar measure
is a state on a $C^\ast$-algebra. In particular, the analogue of the Haar
measure on the deformed $C^\ast$-algebra $\A$ of continuous functions on the
group $SU(2)$ is explicitly known. This Haar functional plays
an important role in the harmonic analysis on the quantum $SU(2)$ group.
For instance, the corepresentations of the $C^\ast$-algebra are
similar to the representations of the Lie group $SU(2)$, and the matrix elements
of the corepresentations can be expressed in terms of the little $q$-Jacobi
polynomials, cf. \cite{\KoorIM}, \cite{\MasuMNNU}, \cite{\VaksS}, and the
orthogonality relations for the little $q$-Jacobi polynomials are
equivalent to the Schur orthogonality relations on the $C^\ast$-algebra
$\A$ involving the Haar functional. This was the start of a fruitful connection
between $q$-special functions and the representation theory of
quantum groups, see e.g. \cite{\KoelAAM}, \cite{\KoorOPTA}, \cite{\Noum}
for more information.

The Haar functional can be restricted to a $C^\ast$-subalgebra of $\A$
generated by a self-adjoint element. It turns out that the Haar
functional restricted to specific examples of such $C^\ast$-subalgebras can
be written as an infinite sum, as an integral with an absolutely continuous
measure or as an integral with an absolutely continuous measure on $[-1,1]$
and a finite number of discrete mass points off $[-1,1]$. These measures
are orthogonality measures for subclasses of the little $q$-Jacobi
polynomials, the big $q$-Jacobi polynomials and the Askey-Wilson polynomials.
These formulas for the Haar functional have been
proved in several cases by Woronowicz \cite{\WoroCMP, App.~A.1} and
Koornwinder \cite{\KoorZSE, thm.~5.3}, see also \cite{\KoorOPTA, thm.~8.4},
by using the representation
theory of the quantum $SU(2)$ group, i.e. the corepresentations of the
$C^\ast$-algebra $\A$ equipped with a suitable comultiplication. See also
Noumi and Mimachi \cite{\NoumMCM, thm.~4.1}.

The proofs by Woronowizc and
Koornwinder use the Schur orthogonality relations for the matrix elements
of irreducible unitary corepresentations of $\A$. They determine a
combination of such matrix elements of one irreducible unitary
corepresentation as an orthogonal polynomial, say $p_n$, in a simple
self-adjoint element, say $\rho$, of the $C^\ast$-algebra $\A$. Here the
degree $n$ of the polynomial $p_n$ is directly related to the spin $l$ of
the irreducible unitary corepresentation. Hence, they conclude that the Haar
functional on the $C^\ast$-subalgebra generated by $\rho$ is given by
an integral with respect to the normalised orthogonality measure for
the polynomials $p_n$. The result by Noumi and Mimachi is closely related to
a limiting case of Koornwinder's result, but the invariant functional lives
on a quantum space on which the quantum $SU(2)$ group acts. Their proof
follows by checking that the moments agree. However, the invariant functional
takes values in a commutative subalgebra of a non-commutative algebra.

It is the purpose of the present paper to prove these results in an
alternative way by only using the $C^\ast$-algebra $\A$. And in particular
we study the spectral properties of the generator of the $C^\ast$-subalgebra
on which the Haar functional is given as a suitable measure. In order to do
so we use the infinite dimensional irreducible representations of the
$C^\ast$-algebra. The infinite dimensional irreducible representations of the
$C^\ast$-algebra $\A$ are parametrised by the unit circle, and the
intersection of the kernels of these representations is trivial. So this set
of representations of $\A$ contains sufficiently many
representations. We determine the spectral properties of the operators
that correspond to the self-adjoint element which generates the
$C^\ast$-subalgebra on which the Haar functional is given by Woronowicz,
Koornwinder and Noumi and Mimachi. An important property of the corresponding
self-adjoint operators is that they can be given as Jacobi matrices, i.e. as
tridiagonal matrices, in a suitable
basis, and hence give rise to orthogonal polynomials. These orthogonal
polynomials can be determined explicitly and can then be used to derive
the explicit form of the Haar functional.

The contents of this paper are as follows. In \S 2 we recall Woronowicz's
quantum $SU(2)$ group and the Haar functional on the corresponding
$C^\ast$-algebra. The spectral theory of the Jacobi matrices is briefly
recalled in \S 3. Woronowicz's \cite{\WoroCMP} expression for the Haar
functional on the algebra of cocentral elements is then proved in \S 4 by
use of the continuous $q$-Hermite polynomials and its Poisson kernel.
In \S 5 we prove the statement of Noumi and Mimachi \cite{\NoumMCM} and
Koornwinder \cite{\KoorZSE} that the Haar functional on a certain
$C^\ast$-subalgebra can be written as a $q$-integral. Here we use
Al-Salam--Carlitz polynomials and $q$-Charlier polynomials. Finally in \S 6
we give a proof of Koornwinder's \cite{\KoorZSE} result that the Haar
functional on certain elements can be written in terms of an Askey-Wilson
integral. Here we use the Al-Salam--Chihara polynomials and the corresponding
Poisson kernel.
It must be noted that the result of \S 5 can be obtained by a formal limit
transition of the result of \S 6, cf. \cite{\KoorZSE, rem.~6.6}, but
we think that the proof in \S 5 is of independent interest, since it is much
simpler. Moreover, the basis of the representation space introduced in \S 5
is essential in \S 6.

To end this introduction we recall some definitions from the theory of
basic hypergeometric series. In this we follow the
excellent book \cite{\GaspR, Ch.~1} by
Gasper and Rahman. We always assume $0<q<1$. The $q$-shifted factorial
is defined by
$$
(a;q)_k=\prod_{i=0}^{k-1}(1-aq^i),\quad
(a_1,\ldots,a_r;q)_k =\prod_{i=1}^r (a_i;q)_k
$$
for $k\in\Zp\cup\{\infty\}$. The $q$-hypergeometric series is defined by
$$
\align
{}_r\vp_s\left( {{a_1,\ldots,a_r}\atop{b_1,\ldots,b_s}};q,z\right) &=
{}_r\vp_s(a_1,\ldots,a_r;b_1,\ldots,b_s;q,z)\\
&= \sum_{k=0}^\infty {{(a_1,\ldots,a_r;q)_k}\over{(q,b_1,\ldots,b_s;q)_k}}
z^k \Bigl( (-1)^k q^{k(k-1)/2}\Bigr)^{s+1-r}.
\endalign
$$
Note that for $a_i=q^{-n}$, $n\in\Zp$, the series terminates, and we find a
polynomial. We also need the $q$-integral;
$$
\int_0^b f(x) \,d_qx = (1-q)b\sum_{k=0}^\infty f(bq^k)q^k,\qquad
\int_a^b f(x)\, d_qx = \int_0^b f(x)\, d_qx - \int_0^a f(x)\, d_qx.
$$
For the very-well-poised ${}_8\vp_7$-series we use the abbreviation, cf.
\cite{\GaspR, Ch.~2},
$$
{}_8W_7(a;b,c,d,e,f;q,z) = {}_8\vp_7 \left(
{{a,q\sqrt{a},-q\sqrt{a},b,c,d,e,f}\atop
{\sqrt{a},-\sqrt{a},qa/b,qa/c,qa/d,qa/e,qa/f}};q,z\right).
\tag\eqname{\vglabbrverywellpoised}
$$

\demo{Acknowledgement} We thank Mizan Rahman and Serge\u\i~Suslov for sending
the preprint \cite{\AskeRS} containing the Poisson kernel for the
Al-Salam--Chihara polynomials. Thanks are also due to Mourad Ismail for
communicating this Poisson kernel.
\enddemo

%%%%%%%%%%%%%%%%%%%%%%%%%%%%%%%%%%%%%%%%%%%%%%%%%%%%%%%%%%%%%%%%%%%%%%
% NEW SECTION %%%%%%%%%%%%%%%%%%%%%%%%%%%%%%%%%%%%%%%%%%%%%%%%%%%%%%%%
%%%%%%%%%%%%%%%%%%%%%%%%%%%%%%%%%%%%%%%%%%%%%%%%%%%%%%%%%%%%%%%%%%%%%%

\head\newsection . The quantum $SU(2)$ group\endhead

We recall in this section Woronowicz's first example of a quantum group,
namely the analogue of the Lie group $SU(2)$, cf. \cite{\WoroRIMS},
\cite{\WoroCMP}. In the general theory of compact matrix quantum groups
Woronowicz has proved the existence of the analogue of a left and right
invariant measure \cite{\WoroCMP}.

We first introduce the $C^{\ast}$-algebra $\A$.
The $C^\ast$-algebra $\A$ is the unital $C^\ast$-algebra
generated by two elements $\a$ and $\g$
subject to the relations
$$
\gathered
\a\g=q\g\a,\qquad \a\g^\ast=q\g^\ast\a, \qquad \g\g^\ast=\g^\ast\g, \\
\a^\ast\a +\g^\ast\g=1=\a\a^\ast + q^2 \g\g^\ast,
\endgathered
\tag\eqname{\vglcomrel}
$$
where $0<q<1$. Here $q$ is a deformation parameter, and for $q=1$ we can
identify $\A$ with the $C^\ast$-algebra of continuous function on
$SU(2)$, where $\a$ and $\g$ are coordinate functions. The group
multiplication is reflected in a $C^\ast$-homomorphism
$\Delta\colon\A\to\A\otimes\A$. (Since $\A$ is a type I
$C^\ast$-algebra, we can take any $C^\ast$-tensor product on the right
hand side.) The Haar functional is then uniquely determined by the
conditions $h(1)=1$ and $(id\otimes h)\Delta(a) = h(a)1 = (h\otimes
id)\Delta(a)$ for all $a\in\A$.

The irreducible representations of the $C^\ast$-algebra $\A$ have been
completely classified, cf. \cite{\VaksS}. Apart from the one-dimensional
representations $\a\mapsto e^{i\theta}$, $\g\mapsto 0$, we only have the
infinite dimensional representations $\p_\phi$, $\phi\in[0,2\p)$,
of $\A$ acting in $\Hi$.
Denote by $\{e_n\mid n\in\Zp\}$ the standard orthonormal basis of $\Hi$,
then $\p_\phi$ is given by
$$
\p_\phi(\a)e_n = \sqrt{1-q^{2n}}e_{n-1},\qquad
\p_\phi(\g) e_n = e^{i\phi}q^n e_n.
\tag\eqname{\vgldefrep}
$$
Here we follow the convention that $e_{-p}=0$ for $p\in\N$.

Woronowicz \cite{\WoroCMP, App.~A.1} has given an explicit formula for
the Haar functional (not using corepresentations) in terms of an infinite
dimensional
faithful representation of $\A$. This can be rewritten in terms of the
irreducible representations of $\A$ as
$$
h(a) = (1-q^2) \sum_{p=0}^\infty {q^{2p}\over{2\p}} \int_0^{2\p} \langle
\p_\phi(a) e_p, e_p\rangle \, d\phi, \quad a\in\A.
$$
Observe that $h\bigl(p(\g^\ast\g)\bigr) = \int_0^1 p(x)\, d_{q^2}x$, for any
continuous $p$ on $\{ q^{2k}\mid k\in\Zp\}$. This
can be considered as a limit case of the results of \S 5, \S 6, cf.
\cite{\KoorZSE, rem.~6.6}.

Introduce the self-adjoint positive
diagonal operator $D\colon\Hi\to\Hi$, $e_p\mapsto q^{2p}e_p$,
then we can rewrite this in a basis independent way;
$$
h(a) = {{(1-q^2)}\over{2\pi}}
\int_0^{2\p} tr\bigl( D\p_\phi(a)\bigr) \, d\phi,\qquad
a\in\A.
\tag\eqname{\vglHaartrace}
$$
The trace operation in \thetag{\vglHaartrace} is well-defined due to the
appearance of $D$. Since $D^{1/2}$ is a Hilbert-Schmidt operator, so
is $\p_\phi(a)D^{1/2}$. The trace of the product of two Hilbert-Schmidt
operators is well-defined. Moreover, $tr\bigl( D\p_\phi(a)\bigr)=
tr\bigl(\p_\phi(a)D\bigr)$ and it is independent of the choice of
the basis. The trace can be estimated by the product of the
Hilbert-Schmidt norms of $D^{1/2}$ and $\p_\phi(a)D^{1/2}$ and then we get
$|tr\bigl(D\p_\phi(a)\bigr)|\leq \parallel a\parallel_{\A} /(1-q^2)$, so that
the function in \thetag{\vglHaartrace} is integrable.
See Dunford and Schwartz \cite{\DunfS, Ch.~XI, \S 6} for more details.

%%%%%%%%%%%%%%%%%%%%%%%%%%%%%%%%%%%%%%%%%%%%%%%%%%%%%%%%%%%%%%%%%%%%%%
% NEW SECTION %%%%%%%%%%%%%%%%%%%%%%%%%%%%%%%%%%%%%%%%%%%%%%%%%%%%%%%%
%%%%%%%%%%%%%%%%%%%%%%%%%%%%%%%%%%%%%%%%%%%%%%%%%%%%%%%%%%%%%%%%%%%%%%

\head\newsection . Spectral theory of Jacobi matrices\endhead

In this section we recall some of the results on the spectral theory
of Jacobi matrices and the relation with orthogonal polynomials. For
more information we refer to Berezanski\u\i\ \cite{\Bere, Ch.~VII, \S 1}
and Dombrowski \cite{\Domb}.

The operator $J$ acting on the standard orthonormal basis
$\{ e_n\mid n\in\Zp\}$ of $\Hi$ by
$$
Je_n =a_{n+1} e_{n+1} + b_n e_n + a_n e_{n-1},\qquad a_n>0,\ b_n\in\R,
\tag\eqname{\vgldefJacobi}
$$
is called a Jacobi matrix. This operator is symmetric, and its deficiency
indices are $(0,0)$ or $(1,1)$. The deficiency indices are $(0,0)$, and thus
$J$ is self-adjoint in $\Hi$, if and only if the corresponding orthonormal
polynomials $p_n(x)$, defined by
$$
xp_n(x) =a_{n+1} p_{n+1}(x) + b_n p_n(x) + a_n p_{n-1}(x),\qquad
p_{-1}(x)=0,\ p_0(x)=1,
\tag\eqname{\vgldeforthonormpols}
$$
correspond to a determined moment problem. In that case we have a unique
measure $dm$ such that
$$
\int_\R p_n(x)p_m(x) \, dm(x) = \d_{n,m}.
$$
All orthogonal polynomials in this paper have compactly supported
orthogonality measures, so that by \cite{\Chih, Ch.~II, thm.~5.6}
the moment problem is determined.

So we assume now that $J$ is self-adjoint.
The point spectrum of $J$ corresponds to the discrete mass points in $dm$.
The spectral decomposition $E$ of $J$ is related to the orthogonality
measure by $\langle E(B)e_n,e_m\rangle = \int_B p_n(x)p_m(x) \, dm(x)$
where $B\subset\R$ is a Borel set.
The operator $J$ can be represented by the multiplication operator on the
Hilbert space $L^2(dm)$ of square integrable functions with respect to the
measure $dm$ as follows. Let $\Lambda\colon\Hi\to L^2(dm)$ be defined
by $e_n\mapsto p_n(x)$. Since the moment problem is determined, the
polynomials are dense in $L^2(dm)$, so that $\Lambda$ can be uniquely
extended to a unitary operator. From \thetag{\vgldefJacobi} and
\thetag{\vgldeforthonormpols} it follows that $\Lambda J v = M \Lambda v$,
$v\in\Hi$,
where $M\colon L^2(dm)\to L^2(dm)$, $Mf(x)=xf(x)$, is the multiplication
operator.
Note that $M$ is bounded if the support of $dm$ is bounded, as is the case
with all orthogonal polynomials considered in this paper.

%%%%%%%%%%%%%%%%%%%%%%%%%%%%%%%%%%%%%%%%%%%%%%%%%%%%%%%%%%%%%%%%%%%%%%
% NEW SECTION %%%%%%%%%%%%%%%%%%%%%%%%%%%%%%%%%%%%%%%%%%%%%%%%%%%%%%%%
%%%%%%%%%%%%%%%%%%%%%%%%%%%%%%%%%%%%%%%%%%%%%%%%%%%%%%%%%%%%%%%%%%%%%%

\head\newsection . The Haar functional on cocentral elements\endhead

Using the characters of the irreducible unitary corepresentations of
the quantum $SU(2)$ group Woronowicz \cite{\WoroCMP, App.~A.1} has
proved an expression for the Haar functional on the $C^\ast$-subalgebra
generated by $\a+\a^\ast$. The set of characters is known as the set
of cocentral elements and it is generated by the element $\a+\a^\ast$.
The purpose of this section is to give a proof of this theorem based on
the spectral analysis of the generator $\a+\a^\ast$. The method used in
this section is also used in section 6 in somewhat greater computational
complexity.

\proclaim{Theorem \theoremname{\thmHaarcocentral}} The Haar functional
on the $C^\ast$-subalgebra generated by the self-adjoint element
$\a+\a^\ast$ is given by the integral
$$
h\bigl(p((\a+\a^\ast)/2)\bigr) = {2\over \p}\int_{-1}^1 p(x)
\sqrt{1-x^2}\, dx
$$
for any continuous function $p\in C\bigl( [-1,1]\bigr)$.
\endproclaim

In order to give an alternative proof of this theorem,
we start with considering
$$
2\p_\phi\big( (\a+\a^\ast)/2\bigr) e_n = \sqrt{1-q^{2n}} e_{n-1} +
\sqrt{1-q^{2n+2}} e_{n+1}.
\tag\eqname{\vglcocentralpi}
$$
So the operator $\p_\phi\big( (\a+\a^\ast)/2\bigr)$ is
represented by a Jacobi matrix with respect to the standard basis of
$\Hi$.

Recall Rogers's continuous $q$-Hermite polynomials $H_n(x| q)$ defined
by
$$
2 x H_n(x| q) = H_{n+1}(x;q) + (1-q^n) H_{n-1}(x| q), \quad
H_{-1}(x| q)=0,\ H_0(x| q)=1.
\tag\eqname{\vglttrqHermite}
$$
The continuous $q$-Hermite polynomials satisfy the orthogonality relations
$$
\int_0^{\p} H_n(\cos\theta| q)H_m(\cos\theta| q)
w(\cos\theta| q) \, d\theta = \d_{n,m} {{2\pi
(q;q)_n}\over{(q;q)_\infty}},
\tag\eqname{\vglorthoqHermite}
$$
with $w(\cos\theta| q) = (e^{2i\theta},e^{-2i\theta};q)_\infty$,
cf. \cite{\AskeI}.
The Poisson kernel for the continuous $q$-Hermite polynomials is given
by, cf. \cite{\AskeI}, \cite{\Bres},
$$
\multline
P_t(\cos\theta,\cos\psi| q) =
\sum_{n=0}^\infty {{H_n(\cos\theta| q)H_n(\cos\psi| q)
t^n}\over{(q;q)_n}} = \\{{(t^2;q)_\infty}\over{(te^{i\theta+i\psi},
te^{i\theta-i\psi},te^{-i\theta+i\psi},te^{-i\theta-i\psi};q)_\infty}}
\endmultline
\tag\eqname{\vglPoissonqHermite}
$$
for $| t| <1$.

Compare \thetag{\vglcocentralpi} with the three-term recurrence relation
\thetag{\vglttrqHermite} to see that \thetag{\vglcocentralpi} is solved
by the orthonormal continuous $q$-Hermite polynomials $H_n(x|
q^2)/\sqrt{(q^2;q^2)_n}$. Using the spectral theory of Jacobi matrices
as in \S 3
we have obtained the spectral decomposition of the self-adjoint operator
$\p_\phi\bigl((\a+\a^\ast)/2\bigr)$, which has spectrum $[-1,1]$.
This link between $\a+\a^\ast$ and the continuous $q$-Hermite
polynomials is already observed in \cite{\KoelAAM, \S 11}. Hence,
for any continuous function $p$ on $[-1,1]$ we have, by the functional
calculus in $C^\ast$-algebras and use of the map $\Lambda$ as in \S 3,
$$
tr\Bigl( D\p_\phi\bigl(p((\a+\a^\ast)/2)\bigr)\Bigr) =
\int_{-1}^1 p(x) P_{q^2}(x,x| q^2) \, dm(x| q^2),
\tag\eqname{\vglcocentralone}
$$
where $dm(x| q^2)=(2\pi)^{-1}(q^2;q^2)_\infty w(x| q^2)
(1-x^2)^{-1/2}dx$ is the corresponding normalised orthogonality measure.

From \thetag{\vglPoissonqHermite} we see that
$$
w(x| q^2) P_{q^2}(x,x| q^2) =
{{4(1-x^2)}\over{(1-q^2) (q^2;q^2)_\infty}},
$$
so that
$$
tr\Bigl( D\p_\phi\bigl(p((\a+\a^\ast)/2)\bigr)\Bigr) =
{2\over{\pi(1-q^2)}} \int_{-1}^1 p(x)\sqrt{1-x^2}\, dx.
$$
Since this is independent of the parameter $\phi$ of the infinite
dimensional representation, we obtain Woronowicz's theorem
\thmref{\thmHaarcocentral} from \thetag{\vglHaartrace}.

%%%%%%%%%%%%%%%%%%%%%%%%%%%%%%%%%%%%%%%%%%%%%%%%%%%%%%%%%%%%%%%%%%%%%%
% NEW SECTION %%%%%%%%%%%%%%%%%%%%%%%%%%%%%%%%%%%%%%%%%%%%%%%%%%%%%%%%
%%%%%%%%%%%%%%%%%%%%%%%%%%%%%%%%%%%%%%%%%%%%%%%%%%%%%%%%%%%%%%%%%%%%%%

\head\newsection . The Haar functional on special spherical
elements\endhead

We start with introducing the self-adjoint element
$$
\rti = iq^\t (\a^\ast\g-\g^\ast\a) - (1-q^{2\t})\g^\ast\g \in \A.
$$
This element is a limiting case of the general spherical element
considered in \S 6, and from that section the notation for this element is
explained. The Haar functional on the $C^\ast$-algebra generated by
$\rti$ can be written as a $q$-integral as proved by Koornwinder
\cite{\KoorZSE, thm.~5.3, rem.~6.6} and Noumi and Mimachi
\cite{\NoumMCM, thm.~4.1} using the corepresentations of the quantum $SU(2)$
group.

\proclaim{Theorem \theoremname{\thmHaaronspecialspherelt}}
The Haar functional
on the $C^\ast$-subalgebra generated by the self-adjoint element
$\rti$ is given by the $q$-integral
$$
h\bigl(p(\rti)\bigr) = {1\over{1+q^{2\t}}}
\int_{-1}^{q^{2\t}} p(x) \, d_{q^2}x
$$
for any continuous function $p$ on $\{ -q^{2k}\mid
k\in\Zp\}\cup\{ q^{2\t+2k}\mid k\in\Zp\}\cup\{ 0\}$.
\endproclaim

To prove this theorem from a spectral analysis of $\p_\phi(\rti)$ we
have to recall the following result.
It is shown in \cite{\KoelAF, prop.~4.1, cor.~4.2}
that $\Hi$ has an orthogonal
basis of eigenvectors of $\p_0(\rti)$. The proof of that proposition
only needs a minor adaption to handle the general case $\p_\phi(\rti)$.
In general we have the following proposition.

\proclaim{Proposition \theoremname{\propeigvectpirti}}
$\Hi$ has an orthogonal basis of
eigenvectors $v_\l^\phi$, where $\l=-q^{2k}$, $k\in\Zp$, and
$\l = q^{2\t + 2k}$,
$k\in\Zp$, for the eigenvalue $\l$ of the self-adjoint operator
$\p_\phi(\rti)$. The squared norm is given by
$$
\align
\langle v_\l^\phi,v_\l^\phi\rangle &= q^{-2k} (q^2;q^2)_k
(-q^{2-2\t};q^2)_k (-q^{2\t};q^2)_\infty,\qquad \l=-q^{2k},\\
\langle v_\l^\phi,v_\l^\phi\rangle &= q^{-2k} (q^2;q^2)_k
(-q^{2+2\t};q^2)_k (-q^{-2\t};q^2)_\infty,\qquad \l=q^{2\t+2k}.
\endalign
$$
Moreover, $v_\l^\phi=\sum_{n=0}^\infty i^ne^{in\phi} p_n(\l)e_n$ with
the polynomial $p_n(\l)$ defined by
$$
\aligned
p_n(\l) &= {{q^{-n\t} q^{\hf n(n-1)}}\over{\sqrt{(q^2;q^2)_n}}}
\, {}_2\vp_1(q^{-2n}, q^{2\t}/\l;0;q;-q^2\l) \\
&= {{(-1)^nq^{n\t} q^{\hf n(n-1)}}\over{\sqrt{(q^2;q^2)_n}}}
\, {}_2\vp_1(q^{-2n}, -1/\l;0;q;q^{2-2\t}\l).
\endaligned
\tag\eqname{\vglexppnl}
$$
\endproclaim

\demo{Remark \theoremname{\rempropeigvectpirti}} The polynomials in $\l$
in \thetag{\vglexppnl} are Al-Salam--Carlitz polynomials $U_n^{(a)}$.
The orthogonality relations obtained from $\langle v_\l^\phi,
v_\m^\phi\rangle = 0$ for $\l\not=\m$ are the orthogonality relations
for the $q$-Charlier polynomials, cf. \cite{\KoelAF, cor.~4.2}. \enddemo

\demo{Remark \theoremname{\remorthogonaldecomp}} The basis described in
proposition \thmref{\propeigvectpirti} induces an orthogonal decomposition
of the representation space $\Hi = V_1^\phi\oplus V_2^\phi$, where
$V_1^\phi$ is the subspace with basis $v^\phi_{-q^{2k}}$, $k\in\Zp$, and
$V_2^\phi$ is the subspace with basis $v^\phi_{q^{2\t+2k}}$, $k\in\Zp$.
\enddemo

In this section we use the basis described in proposition
\thmref{\propeigvectpirti} to calculate the trace. Observe now that for
any function $p$ continuous on the spectrum $\{ -q^{2k}\mid
k\in\Zp\}\cup\{ q^{2\t+2k}\mid k\in\Zp\}\cup\{ 0\}$ we have
$$
tr\Bigl( D\p_\phi\bigl(p(\rti)\bigr)\Bigr) =
\sum_{k=0}^\infty p(-q^{2k}) {{\langle D v_{-q^{2k}}^\phi,
v_{-q^{2k}}^\phi\rangle}\over{\langle v_{-q^{2k}}^\phi,
v_{-q^{2k}}^\phi\rangle}}  +
\sum_{k=0}^\infty p(q^{2\t+2k}){{\langle D v_{q^{2\t+2k}}^\phi,
v_{q^{2\t+2k}}^\phi\rangle}\over{\langle v_{q^{2\t+2k}}^\phi,
v_{q^{2\t+2k}}^\phi\rangle}}.
\tag\eqname{\vgltrsseeen}
$$
So it remains to calculate the matrix coefficients on the diagonal
of the operator $D$ with respect to this basis. We give all the matrix
coefficients in the following lemma.

\proclaim{Lemma \theoremname{\lemmatrixcoeffD}} The matrix coefficients
of the operator $D$ with respect to the orthogonal basis $v_\l^\phi$ are
given by
$$
\align
\langle Dv_{-q^{2k}}^\phi,v_{-q^{2l}}^\phi\rangle &=
(-q^{2\t+2};q^2)_\infty (q^2;q^2)_k (-q^{2-2\t};q^2)_l,
\qquad k\geq l,\\
\langle Dv_{q^{2\t+2k}}^\phi,v_{q^{2\t+2l}}^\phi\rangle &=
(-q^{2-2\t};q^2)_\infty (q^2;q^2)_k (-q^{2+2\t};q^2)_l,
\qquad k\geq l,\\
\langle Dv_{-q^{2k}}^\phi,v_{q^{2\t+2l}}^\phi\rangle &=
(q^2;q^2)_\infty (-q^{2-2\t};q^2)_k (-q^{2+2\t};q^2)_l,
\endalign
$$
and all other cases follow from $D$ being self-adjoint.
\endproclaim

\demo{Proof} The proof is based on calculations involving the
$q$-Charlier polynomials $c_n(x;a;q) =
{}_2\vp_1(q^{-n},x;0;q,-q^{n+1}/a)$. Define a moment functional $\L$
by
$$
\L(p) = \sum_{n=0}^\infty {{q^{2n\t}q^{n(n-1)}}\over{(q^2;q^2)_n}}
p(q^{-2n}),
\tag\eqname{\vgldefL}
$$
for any polynomial $p$.
Note that all moments, i.e. $\L(x^n)$, $n\in\Zp$, exist. The
orthogonality relations $\langle
v^\phi_{-q^{2k}},v^\phi_{-q^{2l}}\rangle =0$ are rewritten as
the orthogonality relations for the $q$-Charlier polynomials;
$$
\L \bigl(c_k(x;q^{2\t};q^2)c_l(x;q^{2\t};q^2)\bigr) =
\d_{k,l} q^{-2k} (q^2,-q^{2-2\t};q^2)_k (-q^{2\t};q^2)_\infty,
$$
cf. \cite{\GaspR, ex.~7.13}, \cite{\KoelAF, cor.~4.2}.

Using proposition \thmref{\propeigvectpirti} and the definition of
the self-adjoint operator $D$ we see that
$$
\langle Dv_{-q^{2k}}^\phi,v_{-q^{2l}}^\phi\rangle =
\L \bigl( {1\over x} c_k(x;q^{2\t};q^2)c_l(x;q^{2\t};q^2)\bigr).
$$
Note that this expression is well-defined, since $0$ is not an element
of the support of the measure representing $\L$. Now use the
orthogonality of the $q$-Charlier polynomials to see that for $k\geq l$
this equals $c_l(0;q^{2\t};q^2) \, \L\bigl(c_k(x;q^{2\t};q^2)/x\bigr)$.
Use the ${}_2\vp_1$-series representation for the $q$-Charlier
polynomials and \thetag{\vgldefL} to evaluate the moment functional
on this particular function. We get
$$
\L\Bigl( {1\over x} c_k(x;q^{2\t};q^2)\Bigr) =
\sum_{l=0}^k {{(q^{-2k};q^2)_l}\over{(q^2;q^2)_l}} (-q^{2k+2+2\t})^l
\sum_{n=l}^\infty {{(q^{-2n};q^2)_l}\over{(q^2;q^2)_n}} q^{2n(\t+1)}
q^{n(n-1)},
$$
where interchanging the summations is allowed as
all sums are absolutely convergent. Replace the summation parameter
$n=m+l$ and use
that $(q^{-2(m+l)};q^2)_l/(q^2;q^2)_{m+l}=(-1)^l q^{-l(l+1+2m)}/(q^2;q^2)_m$.
Now the inner
sum can be summed using ${}_0\vp_0(-;-;q^2,z)=(z;q^2)_\infty$, cf.
\cite{\GaspR, (1.3.16)}. The remaining
finite sum can be summed using the terminating $q$-binomial formula
${}_1\vp_0(q^{-2n};-;q^2,q^{2n}x)=(x;q^2)_n$, cf. \cite{\GaspR, (1.3.14)},
which can also be used to
evaluate the $q$-Charlier polynomial at $x=0$. This proves the first
part of the lemma.
The second part is proved in the same way, but with $\t$ replaced by
$-\t$.

For the last part of the lemma we use the same strategy, but now we have
to use the moment functional $\M$ defined by
$$
\M(p)= \sum_{n=0}^\infty {{(-1)^nq^{n(n-1)}}\over{(q^2;q^2)_n}}
p(q^{-2n}).
$$
for which we have
$\M\bigl( c_k(x;q^{2\t};q^2)c_l(x;q^{-2\t};q^2)\bigr)=0$, cf.
\cite{\KoelAF, cor.~4.2}. \qed\enddemo

Now use lemma \thmref{\lemmatrixcoeffD} and proposition
\thmref{\propeigvectpirti} in \thetag{\vgltrsseeen} to find
$$
tr\Bigl( D\p_\phi\bigl(p(\rti)\bigr)\Bigr) = {1\over{1+q^{2\t}}} \biggl(
\sum_{k=0}^\infty p(-q^{2k}) q^{2k} +
q^{2\t} \sum_{k=0}^\infty p(q^{2\t+2k}) q^{2k}\biggr).
$$
This expression is independent of $\phi$, so that we obtain
from \thetag{\vglHaartrace}
$$
h\bigl(p(\rti)\bigr) = {{1-q^2}\over{1+q^{2\t}}} \biggl(
\sum_{k=0}^\infty p(-q^{2k}) q^{2k} +
q^{2\t} \sum_{k=0}^\infty p(q^{2\t+2k}) q^{2k}\biggr),
$$
which is precisely the statement of theorem
\thmref{\thmHaaronspecialspherelt}.

%%%%%%%%%%%%%%%%%%%%%%%%%%%%%%%%%%%%%%%%%%%%%%%%%%%%%%%%%%%%%%%%%%%%%%
% NEW SECTION %%%%%%%%%%%%%%%%%%%%%%%%%%%%%%%%%%%%%%%%%%%%%%%%%%%%%%%%
%%%%%%%%%%%%%%%%%%%%%%%%%%%%%%%%%%%%%%%%%%%%%%%%%%%%%%%%%%%%%%%%%%%%%%

\head\newsection . The Haar functional on spherical elements\endhead

In this section we give a proof of Koornwinder's theorem expressing the
Haar functional on a $C^\ast$-subalgebra of $\A$ as an Askey-Wilson
integral from the spectral analysis of the generator of the
$C^\ast$-subalgebra.

We first introduce the  self-adjoint element
$$
\multline
\r = \hf\Bigl( \a^2 + (\a^\ast)^2 + q\g^2 + q(\g^\ast)^2 +
iq(q^{-\s} - q^\s)(\a^\ast\g - \g^\ast\a) \\
- iq(q^{-\t} - q^\t)(\g\a - \a^\ast\g^\ast) -
q(q^{-\s} - q^\s)(q^{-\t} - q^\t)\g^\ast\g\Bigr) \in\A.
\endmultline
$$
Note that the element $\rti$, as introduced in section 5, is a limiting case
of this general spherical element, namely
$$
\rti = \lim_{\s\to\infty}2q^{\s +\t -1}\r.
$$
The Haar functional on the $C^\ast$-algebra generated by
$\r$ can be written as an integral with respect to the orthogonality measure
for  Askey-Wilson polynomials, as proved by Koornwinder
\cite{\KoorZSE, thm.~5.3}.

\proclaim{Theorem \theoremname{\thmHaarongeneralspherelt}}
The Haar functional
on the $C^\ast$-subalgebra generated by the self-adjoint element
$\r$ is given by
$$
h(p(\r)) = \int_{\R} p(x)\, dm(x;a,b,c,d| q^2)
\tag\eqname{\vglHaarongeneralspherelt}
$$
for any continuous function $p$ on the spectrum of $\r$, which coincides with
the support of the orthogonality measure in
\thetag{\vglHaarongeneralspherelt}. Here $a=-q^{\s +\t +1}$,
$b=-q^{-\s - \t +1}$, $c=q^{\s -\t +1}$, $d=q^{-\s +\t +1}$ and
$dm(x;a,b,c,d| q^2)$ denotes the normalised Askey-Wilson measure.
\endproclaim

We recall that the normalised Askey-Wilson measure is given by, cf.
Askey and Wilson \cite{\AskeW, thm.~2.1, 2.5},
$$
\int_{\R} p(x)dm(x;a,b,c,d| q) =
{1\over {h_0 2\p}}\int_{0}^{\p}p(\cos\theta)w(\cos\theta)\,d\theta +
{1\over{h_0}} \sum_{k}p(x_k)w_k.
\tag\eqname{\vglHaarongeneralspherelt2}
$$
Here we use the notation $w(\cos\theta)=w(\cos\theta;a,b,c,d|q)$,
$h_0=h_0(a,b,c,d|q)$ and
$$
\aligned
h_0(a,b,c,d|q) & =
{{(abcd;q)_\infty}\over{(q,ab,ac,ad,bc,bd,cd;q)_\infty}}, \\
w(\cos\theta;a,b,c,d|q) & = {{(e^{2i\theta},e^{-2i\theta};q)_\infty} \over
{(ae^{i\theta},ae^{-i\theta}, be^{i\theta},be^{-i\theta},
ce^{i\theta},ce^{-i\theta}, de^{i\theta},de^{-i\theta};q)_\infty}},
\endaligned
\tag\eqname{\vglmeasuregeneral}
$$
and we suppose $a$, $b$, $c$ and $d$ real and such that all pairwise
products are less than $1$.
The sum in \thetag{\vglHaarongeneralspherelt2} is over the points $x_k$
of the form $(eq^k + e^{-1}q^{-k})/2$ with $e$ any of
the parameters $a$, $b$, $c$ or $d$
whose absolute value is larger than one and such that
$|eq^k|>1$, $k\in\Zp$. The corresponding mass
$w_k$ is the residue of $z\mapsto w(\hf(z+z^{-1}))$ at $z = eq^k$ minus
the residue at $z = e^{-1}q^{-k}$. The value of $w_k$ in case $e=a$ is
given in \cite{\AskeW, (2.10)}, but $(1-aq^{2k})/(1-a)$ has to be
replaced by $(1-a^2q^{2k})/(1-a^2)$. Explicitly,
$$
\multline
w_k(a;b,c,d|q) =
{{(a^{-2};q)_\infty}\over{(q,ab,b/a,ac,c/a,ad,d/a;q)_\infty}} \\
\times {{(1-a^2q^{2k})}\over{(1-a^2)}} {{(a^2,ab,ac,ad;q)_k}\over
{(q,aq/b,aq/c,aq/d;q)_k}} \Bigl( {q\over{abcd}}\Bigr)^k,
\endmultline
\tag\eqname{\vgldiscmassAW}
$$
see \cite{\GaspR, (6.6.12)}.

We prove theorem \thmref{\thmHaarongeneralspherelt}
from a spectral analysis of the self-adjoint
operators $\p_\phi(\r)$. We realise $\p_\phi(\r)$ as a Jacobi matrix in
the basis of $\Hi$ introduced in \S 5.

\proclaim{Proposition \theoremname{\propprstJacobimat}} Let $v_\l^\phi$
be the orthogonal basis of $\Hi$ as in proposition
\thmref{\propeigvectpirti}, then
$$
2\p_\phi(\r)v_\l^\phi = qe^{-2i\phi}v_{\l q^2}^\phi + q^{-1}e^{2i\phi}
(1 - q^{-2\t}\l)(1+\l)v_{\l/q^2}^\phi
+ \l q^{1-\t}(q^{-\s}-q^\s)v_\l^\phi
\tag\eqname{\vglrecursiegeneral}
$$
where $\l=-q^{2k}$, $k\in\Zp$, and $\l = q^{2\t + 2k}$, $k\in\Zp$.
\endproclaim

\demo{Proof} We use a factorisation of $\r$ in elements, which are
linear combinations in the generators of the $C^\ast$-algebra $\A$.
Explicitly,
$$
2q^{\t+\s}\r-q^{2\s-1}-q^{2\t+1}=
(\b_{\t+1,\infty}-q^{\s-1}\a_{\t+1,\infty})
(\g_{\t,\infty}+q^\s\d_{\t,\infty}),
\tag\eqname{\vglfactrst}
$$
where
$$
\gathered
\a_{\t,\infty} = q^{1/2}\a+iq^{\t+1/2}\g, \quad
\b_{\t,\infty} = iq^{1/2}\g^\ast + q^{\t-1/2}\a^\ast, \\
\g_{\t,\infty} = -q^{\t+1/2}\a+iq^{1/2}\g, \quad
\d_{\t,\infty} = -iq^{\t+1/2}\g^\ast + q^{-1/2}\a^\ast,
\endgathered
\tag\eqname{\vgldefatinftyetc}
$$
cf. \cite{\Koelpp, prop.~3.3, (2.14), (2.2)}. Of course, \thetag{\vglfactrst}
can also be checked directly from the commutation relations
\thetag{\vglcomrel} in $\A$.
As proved in
\cite{\Koelpp, prop.~3.8}, the operators corresponding to the elements
in \thetag{\vgldefatinftyetc} under the representation $\p_0$
act nicely in the basis $v_\l^0$. The proof immediately generalises to
the basis $v_\l^\phi$. Using the notation $v_\l^\phi(q^\t)=v_\l^\phi$ to
stress the dependence on $q^\t$, we get
$$
\gather
\p_\phi(\a_{\t,\infty})v^\phi_\l(q^\t) =
e^{i\phi}iq^{\hf-\t}(1+\l) v^\phi_{\l/q^2}(q^{\t-1}), \qquad
\p_\phi(\b_{\t,\infty})v^\phi_\l(q^\t) =
e^{-i\phi}iq^\hf v^\phi_\l (q^{\t-1}),\\
\p_\phi(\g_{\t,\infty})v^\phi_\l(q^\t) =
e^{i\phi}iq^\hf (q^{2\t}-\l)v^\phi_{\l}(q^{\t+1}),\qquad
\p_\phi(\d_{\t,\infty})v^\phi_\l(q^\t) =
-e^{-i\phi}iq^{\hf+\t} v^\phi_{\l q^2}(q^{\t+1}).
\endgather
$$
From this result and \thetag{\vglfactrst} the proposition follows.
\qed\enddemo

Proposition \thmref{\propprstJacobimat} implies that
$\p_\phi(\r)$ respects the orthogonal decomposition $\Hi = V_1^\phi
\oplus V_2^\phi$, cf. remark \thmref{\remorthogonaldecomp}.
We denote by $w_m^\phi$, $m\in\Zp$, the orthonormal basis of $V_1^\phi$
obtained by normalising $v_{-q^{2m}}^\phi$, $m\in\Zp$, and by
$u_m^\phi$, $m\in\Zp$, the orthonormal basis of $V_2^\phi$
obtained by normalising $v_{q^{2\t + 2m}}^\phi$, $m\in\Zp$.
Then we get, using proposition \thmref{\propeigvectpirti},
$$
\gathered
2\p_\phi(\r)w_m^\phi = e^{-2i\phi}a_mw_{m+1}^\phi + b_mw_m^\phi + e^{2i\phi}
a_{m-1}w_{m-1}^\phi ,\\
a_m = \sqrt{(1-q^{2m+2})(1+q^{2m+2-2\t})},\qquad b_m = q^{2m+1-\t}
(q^{\s}-q^{-\s})
\endgathered
\tag\eqname{\vglthtrecrstone}
$$
and
$$
\gathered
2\p_\phi(\r)u_m^\phi = e^{-2i\phi}a_mu_{m+1}^\phi + b_mu_m^\phi + e^{2i\phi}
a_{m-1}u_{m-1}^\phi ,\\
a_m = \sqrt{(1-q^{2m+2})(1+q^{2m+2+2\t})},\qquad b_m = q^{2m+1+\t}
(q^{-\s}-q^{\s}).
\endgathered
\tag\eqname{\vglthtrecrsttwo}
$$

In order to solve the corresponding three-term recurrence relations for
the orthonormal polynomials, we recall the Al-Salam--Chihara polynomials
$p_n(x)=p_n(x;a,b|q)$. These polynomials are Askey-Wilson polynomials
\cite{\AskeW} with two parameters set to zero, so the orthogonality
measure for these polynomials follows from
\thetag{\vglHaarongeneralspherelt2}. The Al-Salam--Chihara polynomials are
defined by
$$
p_n(\cos\theta;a,b|q) = a^{-n} (ab;q)_n\, {}_3\vp_2\left(
{{q^{-n},ae^{i\theta},ae^{-i\theta}}\atop{ab,\ 0}};q,q\right).
\tag\eqname{\vgldefAlSalamChiharapols}
$$
The Al-Salam--Chihara polynomials are symmetric in $a$ and $b$
and they satisfy the three-term recurrence relation
$$
2x p_n(x) = p_{n+1}(x) + (a+b)q^n p_n(x) + (1-abq^{n-1})(1-q^n)
p_{n-1}(x).
$$
The Al-Salam--Chihara polynomials are orthogonal with respect to a positive
measure for $ab<1$. Then the orthogonality measure is given by
$dm(x;a,b,0,0|q)$, cf. \thetag{\vglHaarongeneralspherelt2}.
We use the following orthonormal Al-Salam--Chihara polynomials;
$$
h_n(x;s,t| q) = {1\over{\sqrt{(q,-qs^{-2};q)_n}}}
p_n(x;q^{\hf}{t\over s},-q^{\hf}{1\over {st}}| q).
$$
Now the corresponding three-term recurrence relations
\thetag{\vglthtrecrstone}, respectively \thetag{\vglthtrecrsttwo},
are solved by
$e^{2im\phi}h_m(x;q^\t,q^\s| q^2)$, respectively
$e^{2im\phi}h_m(x;q^{-\t},q^{-\s}| q^2)$.
So we see that $\p_\phi(\r)$ preserves the orthogonal decomposition
$\Hi=V_1^\phi\oplus V_2^\phi$ of the representation space and that this
operator is represented by a Jacobi matrix on each of the components.

The operator $D\colon\Hi\to\Hi$ introduced in \S 2 does not preserve the
orthogonal decomposition $\Hi=V_1^\phi\oplus V_2^\phi$ as follows from
lemma \thmref{\lemmatrixcoeffD}. Let
$$
D=\pmatrix D^\phi_{1,1}&D^\phi_{1,2}\\ D^\phi_{2,1}&D^\phi_{2,2}\endpmatrix
$$
be the corresponding decomposition of $D$, then
$$
tr\bigl( D\p_\phi(p(\r))\bigr) =
tr_{V_1^\phi}\bigl( D_{11}^\phi\p_\phi(p(\r))\bigr) +
tr_{V_2^\phi}\bigl( D_{22}^\phi\p_\phi(p(\r))\bigr),
\tag\eqname{\spoor}
$$
since $\p_\phi(p(\r))$ does preserve the orthogonal decomposition.
Moreover, as in \S 3, we define
$$
\align
&\Lambda_1\colon V_1^\phi \to L^2(dm_1), \qquad w^\phi_m\mapsto
e^{2im\phi}h_m(x;q^\t,q^\s| q^2),\\
&\Lambda_2\colon V_2^\phi \to L^2(dm_2), \qquad u^\phi_m\mapsto
e^{2im\phi}h_m(x;q^{-\t},q^{-\s}| q^2),
\endalign
$$
where
$$
dm_1(x)  = dm(x;q^{1+\s-\t},-q^{1-\s-\t},0,0| q^2), \
dm_2(x)  = dm(x;q^{1-\s+\t},-q^{1+\s+\t},0,0| q^2)
$$
are the corresponding normalised orthogonality measures.
Then, by using the spectral theory of Jacobi matrices as described
in \S 3, we get
$$
\multline
tr_{V_1^\phi}\bigl( D_{11}^\phi\p_\phi(p(\r))\bigr) =  \\
\int_{\R}p(x) \sum_{n=0}^{\infty} \sum_{m=0}^{\infty}
\langle  D w_n^\phi,w_m^\phi\rangle h_n(x;q^\t,q^\s| q^2)
h_m(x;q^\t,q^\s| q^2)e^{2i(m-n)\phi}\,dm_1(x)
\endmultline
\tag\eqname{\vglspooreen}
$$
and similarly
$$
\multline
tr_{V_2^\phi}\bigl( D_{22}^\phi\p_\phi(p(\r))\bigr) =   \\
\int_{\R}p(x) \sum_{n=0}^{\infty} \sum_{m=0}^{\infty}
\langle  D u_n^\phi,u_m^\phi\rangle h_n(x;q^{-\t},q^{-\s}| q^2)
h_m(x;q^{-\t},q^{-\s}| q^2)e^{2i(m-n)\phi}\,dm_2(x).
\endmultline
\tag\eqname{\vglspoortwee}
$$
The double sum in both \thetag{\vglspooreen} and \thetag{\vglspoortwee}
is absolutely convergent, uniformly in $\phi$ and uniformly in $x$ on the
support of the corresponding orthogonality measure. To see this we observe
that proposition \thmref{\propeigvectpirti} and lemma
\thmref{\lemmatrixcoeffD} imply
$$
|\langle Dw_n^\phi,w_m^\phi\rangle| \leq
q^{n+m} (-q^2,-q^{2-2\t};q^2)_\infty/(q^2;q^2)_\infty,
$$
so that for some constant $C$
$$
\Bigl\vert \sum_{n,m=0}^\infty
\langle Dw_n^\phi,w_m^\phi\rangle \bigl(h_nh_m\bigr)(x;q^\t,q^\s|q^2)
e^{2i(m-n)\phi} \Bigr\vert \leq
C \left( \sum_{n=0}^\infty q^n |h_n(x;q^\t,q^\s|q^2)|\right)^2.
$$
Now we use the asymptotic behaviour of the Al-Salam--Chihara polynomials
on $[-1,1]$, cf. Askey and Ismail \cite{\AskeIMAMS, \S 3.1}, and, if
$(aq^k+a^{-1}q^{-k})/2$ is a discrete mass point of the orthogonality
measure of the Al-Salam--Chihara polynomials we have
$p_n((aq^k+a^{-1}q^{-k})/2;a,b|q)\sim a^{-n}q^{-nk}$ as $n\to\infty$, cf.
\thetag{\vgldefAlSalamChiharapols},
to find our assertion. The other double sum is treated analogously.

Now that we have established \thetag{\spoor}, \thetag{\vglspooreen} and
\thetag{\vglspoortwee}, we have an explicit expression for the Haar functional
on the $C^\ast$-subalgebra of $\A$ generated by $\r$. Integrate
\thetag{\vglspooreen} and \thetag{\vglspoortwee} with respect to $\phi$ over
$[0,2\pi]$, cf. \thetag{\vglHaartrace},
and interchange summations, which is jusitified by the previous
remark. Since the inner products in \thetag{\vglspooreen} and
\thetag{\vglspoortwee} are independent of $\phi$ by lemma
\thmref{\lemmatrixcoeffD}, the integration over $\phi$
reduces the double sum to a single sum. Now use
\thetag{\vglHaartrace}, lemma \thmref{\lemmatrixcoeffD} and proposition
\thmref{\propeigvectpirti} to prove the following proposition.

\proclaim{Proposition \theoremname{\propHaarandPoisson}}
For any continuous function $p$ on the spectrum of $\r\in\A$, i.e. the
union of the supports of the measures $dm_1$ and $dm_2$,
we have the following expression for the Haar functional
$$
\multline
h\bigl(p(\r)\bigr) = {{1-q^2}\over {1+q^{2\t}}}
\int_\R p(x) P_{q^2}(x,x;q^{1+\s-\t},-q^{1-\s-\t}\vert q^2) \, dm_1(x)
\\ + {{1-q^2}\over {1+q^{-2\t}}}
\int_\R p(x) P_{q^2}(x,x;q^{1-\s+\t},-q^{1+\s+\t}\vert q^2)\, dm_2(x),
\endmultline
\tag\eqname{\vglHaarinterm1}
$$
where
$$
P_t(x,y;a,b\vert q) =
\sum_{k=0}^\infty t^k {{p_k(x;a,b\vert q)p_k(y;a,b
\vert q)} \over {(q,ab;q)_k}}
\tag\eqname{\vgkdefPoissonkernel}
$$
is the Poisson kernel for the Al-Salam--Chihara polynomials
\thetag{\vgldefAlSalamChiharapols}.
\endproclaim

\demo{Remark \theoremname{\remonconvPKAlSCp}} For all values needed the
Poisson kernel in \thetag{\vglHaarinterm1}
is absolutely convergent for $t=q^2$ by the asymptotic
behaviour of the Al-Salam--Chihara polynomials, cf. the remarks following
\thetag{\vglspooreen} and \thetag{\vglspoortwee}.
\enddemo

In order to tie \thetag{\vglHaarinterm1} to
theorem \thmref{\thmHaarongeneralspherelt}
we have to use the explicit expression for the
Poisson kernel of the Al-Salam--Chihara polynomials given by Askey, Rahman and
Suslov \cite{\AskeRS, (14.8)}. It is given in terms of a very-well-poised
${}_8\vp_7$-series, cf. the notation \thetag{\vglabbrverywellpoised},
$$
\multline
P_t(\cos\theta,\cos\psi;a,b\vert q) =
{{(ate^{i\theta},ate^{-i\theta},bte^{i\psi},bte^{-i\psi},t;q)_\infty }\over
{(te^{i\theta +i\psi},te^{i\theta -i\psi},te^{i\psi -i\theta},
te^{-i\psi -i\theta} ,abt;q)_\infty}} \\
\times {}_8W_7({{abt}\over q};t,be^{i\theta},be^{-i\theta},ae^{i\psi}
,ae^{-i\psi};q,t),
\endmultline
\tag\eqname{\PoissonkernelAlSalamChihara}
$$
where we also used the transformation formula
\cite{\GaspR, (2.10.1)}. Askey, Rahman and Suslov \cite{\AskeRS} are not very
specific about the validity of \thetag{\PoissonkernelAlSalamChihara}, but
from \cite{\AskeRS, \S 1} we may deduce that it is valid for $|a|<1$, $|b|<1$
and $|t|<1$. We first observe that \thetag{\PoissonkernelAlSalamChihara}
also holds for $ab<1$ and $|t|<1$, for which
\thetag{\vgkdefPoissonkernel} is absolutely convergent. To see this we
show that the simple poles of the infinite product on the right hand side
\thetag{\PoissonkernelAlSalamChihara} at $1/t=abq^l$, $l\in\Zp$, are
cancelled by a zero of the very-well-poised ${}_8\vp_7$-series at
$1/t=abq^l$, i.e.
$$
{}_8W_7(q^{-l-1};{{q^{-l}}\over{ab}},be^{i\theta},be^{-i\theta},ae^{i\psi}
,ae^{-i\psi};q,{{q^{-l}}\over{ab}}) = 0, \qquad l\in\Zp.
$$
This follows by writing the very-well-poised ${}_8\vp_7$-series as a sum of
two balanced ${}_4\vp_3$-series by \cite{\GaspR, (2.10.10)}, in which case
both ${}_4\vp_3$-series have a factor $(q^{-l};q)_\infty=0$ in front.

We also need the Poisson kernel evaluated in possible discrete mass points
of the corresponding orthogonality measure. For $x$ in a bounded set,
containing the support of the orthogonality measure, we see that for
$|t|$ small enough $P_t(x,x;a,b|q)$ defined by \thetag{\vgkdefPoissonkernel}
is absolutely convergent and uniformly in $x$, and by analytic continuation
in $x$ we see that \thetag{\PoissonkernelAlSalamChihara} is valid for $|t|$
small enough.
Now assume that $|a|>1$ and that $x_k=(aq^k+a^{-1}q^{-k})/2$ is a discrete
mass point of the orthogonality measure of the Al-Salam--Chihara polynomials,
then the radius of convergence of $P_t(x_k,x_k;a,b|q)$
in \thetag{\vgkdefPoissonkernel}
is $a^2q^{2k}$, so the radius of convergence is greater than $1$, cf.
\thetag{\vglHaarongeneralspherelt2}. For this choice of the arguments
the right hand side of \thetag{\PoissonkernelAlSalamChihara} can be expressed
in terms of a terminating very-well-poised ${}_8\vp_7$-series;
$$
(a^2q^kt,tq^{-k};q)_k {{(abtq^k,btq^{-k}/a;q)_\infty}\over
{(abt,ta^{-2}q^{-2k};q)_\infty}} \,
{}_8W_7 ({{abt}\over q};t, abq^k,bq^{-k}/a,q^{-k},a^2q^k;q,t).
\tag\eqname{\vglPkatdiscretemass}
$$
This expression has no poles in the disc $|t|<a^2q^{2k}$, and coincides
with the Poisson kernel \thetag{\vgkdefPoissonkernel}. So we have shown that
\thetag{\PoissonkernelAlSalamChihara} for $t=q$ with
$x=\cos\theta=\cos\psi$ is valid for $ab<1$ and $x$ in the support of the
corresponding orthogonality measure.

After these considerations on the Poisson kernel for the Al-Salam--Chihara
polynomials, we can use Bailey's summation formula, cf.
\cite{\GaspR, (2.11.7)}, in the following form;
$$
\multline
{1\over {(b/a;q)_\infty}} {}_8W_7(a;b,c,d,e,f ;q,q) +
{1\over {(a/b ;q)_\infty}} \\
\times {{(a q, c ,d ,e ,f ,bq/ c , bq/d ,bq/e, bq/f;q)_\infty}\over
{(aq/c, aq/d, aq/e, aq/f, bc/a, bd/a, be/a, bf/a,
b^2q/a;q)_\infty}} \,  {}_8W_7({b^{ 2}\over a };b ,
{{b c }\over a },{{b d }\over a }, {{b e }\over a },{{b f }\over a } ;q,q) \\
= {{(aq, aq/(cd), aq/(ce), aq/(cf), aq/(de), aq/(df), aq/(ef);q)_\infty}\over
{(aq/c, aq/d, aq/e, aq/f, bc/a, bd/a, be/a, bf/a; q)_\infty}}.
\endmultline
$$
Use Bailey's formula
with $q$ replaced by $q^2$ and parameters $a = -q^{2-2\t}$, $b = q^2$, $c =
- q^{1-\s-\t}e^{i\theta}$, $d = - q^{1-\s-\t}e^{-i\theta}$,
$e = q^{1+\s-\t}e^{i\theta}$ and $f = q^{1+\s-\t}e^{-i\theta}$, and
multiply the resulting identity by
$$
{{(1-q^2)(1-e^{2i\theta})(1-e^{-2i\theta})\,
(q^2,-q^{2+2\t},-q^{2-2\t};q^2)_\infty}\over
{(1-q^{1+\s-\t}e^{i\theta})(1-q^{1+\s-\t}e^{-i\theta})
(1+q^{1-\s-\t}e^{i\theta})(1+q^{1-\s-\t}e^{-i\theta})}}
$$
to find, using the notation of \thetag{\vglmeasuregeneral},
$$
\multline
{{(1-q^2) P_{q^2}(x,x;q^{\s -\t +1},-q^{-\s - \t +1}\vert q^2)
w(x;q^{\s-\t+1},-q^{1-\s-\t},0,0|q^2)}\over{
(1+q^{2\t})h_0(q^{\s-\t+1},-q^{1-\s-\t},0,0|q^2)}} + \\
{{(1-q^2) P_{q^2}(x,x;q^{-\s +\t +1},-q^{\s +\t +1}\vert q^2)
w(x;q^{\t-\s+1},-q^{1+\s+\t},0,0|q^2)}\over
{(1-q^{-2\t})h_0(q^{\t-\s+1},-q^{1+\s+\t},0,0|q^2)}} = \\
{{w(x;q^{\t-\s+1},-q^{1+\s+\t},q^{\s-\t+1},-q^{1-\s-\t}|q^2)}
\over{h_0(q^{\t-\s+1},-q^{1+\s+\t},q^{\s-\t+1},-q^{1-\s-\t}|q^2)}}
\endmultline
\tag\eqname{\vglBaileynaarHaar}
$$
for $x\in [-1,1]$, the absolutely continuous part of $dm_1$ and $dm_2$.
This proves the absolutely continuous part of the Askey-Wilson integral
in theorem \thmref{\thmHaarongeneralspherelt}.

In order to deal with the possible discrete mass points in the orthogonality
measures on the right hand side of \thetag{\vglHaarinterm1}, we first observe
that the discrete mass points of $dm_1$ do not occur as discrete mass points
of $dm_2$ and vice versa. Furthermore, note that for $t=q$ the
expression in \thetag{\vglPkatdiscretemass} reduces to
$$
P_k(a;b|q) = {{(abq^k,bq/a;q)_\infty}\over{(ab,a^{-2}q^{1-2k};q)_\infty}}
(q^{-k},bq^{-k}/a,a^2q^k;q)_k q^k.
$$
This can be seen directly from \thetag{\vglPkatdiscretemass}, since only the
last term in the ${}_8\vp_7$-series survives, or by applying Jackson's
summation formula \cite{\GaspR, (2.6.2)}. Now a straightforward calculation
using this formula and the explicit values for the
weights given in \thetag{\vgldiscmassAW} proves
$$
{{1-q}\over{1-q/(ab)}} {{P_k(a;b|q) w_k(a;b,0,0|q)}\over{h_0(a,b,0,0|q)}} =
{{w_k(a;b,q/a,q/b|q)}\over{h_0(a,b,q/a,q/b|q)}}.
$$
This proves that the discrete mass
points in \thetag{\vglHaarinterm1} lead to the discrete mass points in theorem
\thmref{\thmHaarongeneralspherelt}. This proves theorem
\thmref{\thmHaarongeneralspherelt} from the spectral analysis of
$\p_\phi(\r)$.

\demo{Remark \theoremname{\remonresidues}}
It is not allowed to take residues in \thetag{\vglBaileynaarHaar} to prove
the statement concerning the discrete mass points. For this we have to know
that \thetag{\vglBaileynaarHaar} also holds in a neighbourhood of the
discrete mass point $x_k$, but the explicit expression for the Poisson kernel
in \thetag{\PoissonkernelAlSalamChihara} leading  to
\thetag{\vglBaileynaarHaar} may fail to hold.
\enddemo

%%%%%%%%%%%%%%%%%References%%%%%%%%%%%%%%%%%%%%%%%%%%%%%%%%%%%%%%%%%%%%%%

\enddocument